\documentclass[a4paper,10pt]{article}
\usepackage{amssymb,amsmath}
\newcommand{\twosum}[2]{\sum_{\begin{array}{c} {\scriptstyle #1}\\
{\scriptstyle #2} \end{array}}}
\newcommand{\nd}{\, | \hspace{-1.1mm}/\,}
\newcommand{\vp}{\varpi}
\newcommand{\ep}{\varepsilon}

\newcommand{\N}{\mathbb{N}}

\newcommand{\Z}{\mathbb{Z}}

\newcommand{\cl}[1]{{\cal #1}}
\newtheorem{theorem}{Theorem}
\newtheorem{lemma}{Lemma}

\begin{document}
\title{Integers Represented as a Sum of Primes and Powers of Two}
\author{D.R. Heath-Brown and J.-C. Puchta\\
Mathematical Institute, Oxford}
\date{}
\maketitle

\section{Introduction}

It was shown by Linnik [10] that there is an absolute constant $K$ such
that every sufficiently large even integer can be written as a sum of
two primes and at most $K$ powers of two.   This is a remarkably strong
approximation to the Goldbach Conjecture.  It gives us a
very explicit set $\cl{K}(x)$ of integers $n\le x$ of cardinality only
$O((\log x)^K)$, such that every sufficiently large even integer 
$N\le x$ can be written as $N=p+p'+n$, with  $p,p'$ prime and
$n\in\cl{K}(x)$.  In contrast, if one tries to arrange such a
representation using an interval in place of the set $\cl{K}(x)$, all 
known results would require $\cl{K}(x)$ to
have cardinality at least a positive power of $x$.

Linnik did not establish an explicit value for the number $K$ of powers
of 2 that would be necessary in his result.  However, such a value has
been computed by Liu, Liu and Wang [12], who found that $K=54000$ is 
acceptable.  This result
was subsequently improved, firstly by Li [8] who obtained $K=25000$,
then by Wang [18], who found that $K=2250$ is acceptable, and finally
by Li [9] who gave the value $K=1906$. 
One can do better if one assumes the Generalized Riemann Hypothesis, and
Liu, Liu and Wang [13] showed that $K=200$ is then admissible.

The object of this paper is to give a rather different approach to this
problem, which leads to dramatically improved bounds on the number of
powers of 2 that are required for Linnik's theorem.

\begin{theorem}
Every sufficiently large even integer is a sum of two primes and exactly
13 powers of 2.
\end{theorem}
\begin{theorem}
Assuming the Generalized Riemann Hypothesis, every sufficiently 
large even integer is a sum of two primes and exactly
7 powers of 2.
\end{theorem}

We understand that Ruzsa and Pintz have, in work in preparation, 
given an independent proof of Theorem 2, and have established a
version of Theorem 1 requiring only 8 powers of 2.  Indeed, 
already in 2000, Pintz had announced the values $K=12$ unconditionally, and
$K=10$ on the Generalized Riemann Hypothesis.  Although we have not seen
an account of this work, we understand that that our approach is 
different in a number of respects.  

We should also report that Elsholtz, in unpublished work, has shown that 
one can obtain $K=12$ in Theorem 1, by a variant of our method. 
He does this by improving our constant
$2.7895$ in (25) to $2.96169$, by using $D=21$, 
and replacing our estimate (41) for $C_2$ 
by $C_2\le 1.992$.

Previous workers have based their line of attack on a proof of Linnik's
theorem due to Gallagher [3].  Let $\vp$ be a small positive
constant.  Set
\begin{equation}
S(\alpha)=\sum_{\vp N<p\le N}e(\alpha p),
\end{equation}
where $e(x):=\exp(2\pi ix)$,
and 
\[T(\alpha)=\sum_{1\le\nu\le L}e(\alpha 2^{\nu}),\;\;\;
L=[\frac{\log N/2K}{\log 2}].\]

As in earlier proofs of Linnik's Theorem we shall use estimates for ${\rm
meas}(\cl{A}_{\lambda})$, where
\[\cl{A}_{\lambda}=\{\alpha\in[0,1]: |T(\alpha)|\ge\lambda L\}.\]
In \S 7 we shall bound ${\rm meas}(\cl{A}_{\lambda})$ by a new
method, suggested to us by Professor Keith Ball.  This provides
the following estimates.
\begin{lemma}
We have
\[{\rm meas}(\cl{A}_{\lambda})\ll N^{-E(\lambda)}\]
with $E(0.722428)>1/2$ and $E(0.863665)>109/154$.
\end{lemma}
We are extremely grateful to
Professor Ball for suggesting his alternative approach to us.
An earlier version of this paper used a completely different technique
to bound $E(\lambda)$ and showed that one can take
\[E(\lambda)\geq 0.822\lambda^2 +o(1)\]
as $N\rightarrow\infty$.  This sufficed to establish Theorems 1 and 2
with 24 and 9 powers of 2 respectively.  

For comparison with Lemma 1, the best bound for $E(\lambda)$
in the literature is due to Liu, Liu and Wang [11; Lemma 3], and states
that 
\[E(1-\eta)\le 1-F(\frac{2+\sqrt{2}}{4}\eta)-F(1-\frac{2+\sqrt{2}}{4}\eta)
+o(1)\]
for $\eta<(7e)^{-1}$, where $F(x)=x(\log x)/(\log 2)$. 

The estimate provided by Lemma 1 will be injected into the circle
method, where it will be crucial in bounding the minor arc contribution.
On the major arcs we shall improve on Gallagher's analysis so as to show
that hypothetical zeros close to $\sigma=1$ play no r\^{o}le.  Thus,
in contrast to previous workers, we
will have no need for explicit numerical zero-free regions for
$L$-functions.  Naturally this produces a considerable simplification in
the computational aspects of our work.  Thus it is almost entirely the 
values of the constants in Lemma 1 which determine the number of
powers of 2 appearing in Theorems 1 and 2.

The paper naturally divides into two parts, one of which 
involves the circle method and zeros of $L$-functions, and the other
of which is devoted to the proof of Lemma 1.  We begin
with the former.

One remark about notation is in order.  At various stages in the
proof, numerical upper bounds on $\vp$ will
be required.  Since we shall always take $\vp$ to be sufficiently
small, we shall assume that any such bound is satisfied.  Moreover, 
since $\vp$
is to be thought of as fixed, we will allow the implied constants in
the $O(\ldots)$ and $\ll$ notations to
depend on $\vp$.

\section{The Major Arcs}

We shall follow the method of Gallagher [3; \S 1] closely.  We choose a
parameter $P$ in the range $1\le P\le N^{2/5}$ and define the major arcs
${\mathfrak M}$ as the set of $\alpha\in[0,1]$ for which there exist $a\in\Z$
and $q\in \N$ such that $q\le P$ and
\[|\alpha-\frac{a}{q}|\le\frac{P}{qN}.\]
If $\chi$ is a character to modulus $q$, we write
\[c_n(\chi)=\sum_{a=1}^{q}\chi(a)e(\frac{an}{q})\]
and
\[\tau(\chi)=\sum_{a=1}^{q}\chi(a)e(\frac{a}{q}).\]
Moreover we put
\[A(\chi,\beta)=\sum_{\vp N<p\le N}\chi(p)e(\beta p)\]
and
\[I_{n,s}(\chi,\chi')=\int_{-P/sN}^{P/sN}A(\chi,\beta)
A(\chi',\beta)e(-\beta n)d\beta.\]
If $\chi$ is a character to a modulus $r|q$ we also write $\chi_q$ for the
induced character modulo $q$, and if $\chi,\chi'$ are characters to
moduli $r$ and $r'$ respectively, we set
\[J_n(\chi,\chi')=\twosum{q\le P}{[r,r']|q}\frac{1}{\phi(q)^2}
c_n(\chi_q\chi'_q)\tau(\overline{\chi_q})\tau(\overline{\chi'_q})
I_{n,q}(\chi,\chi').\]
Then, by a trivial variant of
the argument leading to Gallagher [3; (3)], we find that 
\begin{equation}
\int_{{\mathfrak M}}S(\alpha)^2 e(-\alpha n)d\alpha=
\sum_{\chi,\chi'}J_n(\chi,\chi')+O(P^{5/2}),
\end{equation}
for any integer $n$, the sum being over primitive characters
$\chi,\chi'$ to moduli $r,r'$ for which $[r,r']\le P$.  
In what follows we shall take $1\le
n\le N$.

To estimate the contribution from a particular pair of characters
$\chi,\chi'$ we put
\[A_q(\chi)=\{\int_{-P/qN}^{P/qN}|A(\chi,\beta)|^2 d\beta\}^{1/2}\]
and
\[C_n(\chi,\chi')=\twosum{q\le P}{[r,r']|q}\frac{1}{\phi(q)^2}
|c_n(\chi_q\chi'_q)\tau(\overline{\chi_q})\tau(\overline{\chi'_q})|.\]
Note that what Gallagher calls $||A(\chi)||$ is our $A_1(\chi)$.  We
have $A_q(\chi)\le A_m(\chi)$ whenever $m\le q$.  Then, as in Gallagher
[3; (4)] we find
\begin{equation}
|J_n(\chi,\chi')|\le C_n(\chi,\chi')A_{[r,r']}(\chi)A_{[r,r']}(\chi').
\end{equation}
It is in bounding $C_n(\chi,\chi')$ that
there is a loss in Gallagher's argument.  Let $r''$ be the conductor
of $\chi\chi'$, and write $m=[r,r']$.  moreover, for any positive
integers $a$ and $n$ we write
\[a_n=\frac{a}{(a,n)}.\]
Then Gallagher shows that
\[C_n(\chi,\chi')\le (rr'r'')^{1/2}\sum_{q\le
P,\,m|q}(\phi(q)\phi(q_n))^{-1},\]
where $q/m$ is square-free and coprime to $m$.  Moreover we have $r''|m_n$.
It follows that
\[C_n(\chi,\chi')\le \frac{(rr'r'')^{1/2}}{\phi(m)\phi(m_n)}
\sum_{(s,m)=1}\mu^2 (s)/\phi(s)\phi(s_n).\]
The sum on the right is
\[\prod_{p\nd mn}(1+\frac{1}{(p-1)^2})
\prod_{p|n, p\nd m}(1+\frac{1}{(p-1)})\ll
\prod_{p|n, p\nd m}\frac{p}{(p-1)},\]
and
\[\frac{m}{\phi(m)}\prod_{p|n, p\nd m}\frac{p}{(p-1)}\le
\frac{n}{\phi(n)}\frac{m_n}{\phi(m_n)}.\]
We therefore deduce that
\[C_n(\chi,\chi')\ll \frac{(rr'r'')^{1/2}}{m}\frac{m_n}{\phi^2 (m_n)}
\frac{n}{\phi(n)}.\]
Now if $p^e||r$ and $p^f||r'$, then $p^{|e-f|}|r''$, since $r''$ is
the conductor of $\chi\chi'$.  (Here the notation $p^e||r$ means, as 
usual, that $p^e|r$ and $p^{e+1}\nd r$.)  We
therefore set 
\begin{equation}
h=(r,r')\;\;\;\mbox{and}\;\;\; r=hs,\; r'=hs', 
\end{equation}
so that $ss'|r''$ and $m=hss'$.  Since
\[\frac{m_n}{\phi^2(m_n)}\ll m_n ^{\vp-1}\]
we therefore have
\[\frac{(rr'r'')^{1/2}}{m}\frac{m_n}{\phi^2 (m_n)}\ll
(ss')^{-1/2}{r''}^{1/2}m_n^{\vp-1}.\]
Now, using the bounds $r''\le m_n$ and $ss'\le r''$, 
we find that
\begin{eqnarray*}
\frac{(rr'r'')^{1/2}}{m}\frac{m_n}{\phi^2 (m_n)}&\ll&
(ss')^{-1/2}{r''}^{1/2}{r''}^{\vp-1}\\
&=&(ss')^{-1/2}{r''}^{\vp-1/2}\\
&\ll & (ss')^{\vp-1}.
\end{eqnarray*}
Alternatively, using only the fact that $m_n\ge r''$, we have
\begin{eqnarray*}
\frac{(rr'r'')^{1/2}}{m}\frac{m_n}{\phi^2 (m_n)}&\ll & 
(ss')^{-1/2}m_n^{1/2}m_n^{\vp-1}\\
&\ll & m_n^{\vp-1/2}.
\end{eqnarray*}

These estimates produce
\[C_n(\chi,\chi')\ll \min\{(ss')^{\vp-1}\,,\,m_n^{\vp-1/2}\}
\frac{n}{\phi(n)}.\]
On combining this with the bounds (2) and (3) we deduce the following
result.
\begin{lemma}
Suppose that $P\le N^{2/5-\vp}$.  Then
\[\int_{{\mathfrak M}}S(\alpha)^2 e(-\alpha n)d\alpha=J_n(1,1)+
O(\frac{n}{\phi(n)}S_n)+O(N^{1-\vp}),\]
where
\[S_n=\sum_{\chi,\chi'}A_{[r,r']}(\chi)
A_{[r,r']}(\chi')\min\{(ss')^{\vp-1}\,,\,m_n^{-1/3}\},\]
the sum being over primitive characters, not both principal, of moduli
$r,r'$, with $[r,r']\le P$.  
\end{lemma}

We have next to consider $A_m(\chi)$.  According to the argument of
Montgomery and Vaughan [15; \S 7] we have
\[A_m(\chi)\ll N^{1/2}\max_{\vp N<x\le N}\max_{0<h\le x}(h+mN/P)^{-1}
|\sum_{x}^{x+h}\chi(p)|.\]
Note that we have firstly taken account of the restriction in (1) to
primes $p>\vp N$, and secondly replaced $(h+N/P)^{-1}$ as it occurs in
Montgomery and Vaughan, by the smaller quantity $(h+mN/P)^{-1}$.  The
argument of [15; \S 7] clearly allows this.

By partial summation we have
\[\sum_{x}^{x+h}\chi(p)\ll 
(\log x)^{-1}\max_{0<j\le h}\sum_{x}^{x+j}\chi(p)\log p.\]
Moreover, a standard application of the `explicit formula' for $\psi(x,\chi)$
produces the estimate
\[\sum_{x}^{x+j}\chi(p)\log p\ll N^{1/2+3\vp}(\log N)^2+
\sum_{\rho}|\frac{(x+j)^{\rho}}{\rho}-\frac{x^{\rho}}{\rho}|,\]
where the sum over $\rho$ is for zeros of $L(s,\chi)$ in the
region
\[\beta\ge\frac{1}{2}+3\vp,\;\;\;|\gamma|\le N.\]
When $\chi$ is the trivial character we shall include the pole $\rho=1$ 
amongst the `zeros'.
Since $j\le h$ and
\[\frac{(x+j)^{\rho}}{\rho}-\frac{x^{\rho}}{\rho}\ll 
\min\{jN^{\beta-1}\,,\,N^{\beta}|\gamma|^{-1}\},\]
we find that
\[A_m(\chi)\ll  \frac{P}{m}N^{4\vp}+\frac{N^{1/2}}{\log N}
\{\max_{0<h\le N}(h+mN/P)^{-1}
\sum_{\rho}N^{\beta-1}\min\{h\,,\,N|\gamma|^{-1}\}.\]
However we have
\[\min\{\frac{h}{h+H}\,,\,\frac{A}{h+H}\}\le\min\{1\,,\,\frac{A}{H}\}\]
whenever $h,H,A>0$.  Applying this with $H=mN/P$ and
$A=N|\gamma|^{-1}$, we deduce that
\begin{equation}
A_m(\chi)\ll  \frac{P}{m}N^{4\vp}+\frac{N^{1/2}}{\log N}
\sum_{\rho}N^{\beta-1}\min\{1\,,\,Pm^{-1}|\gamma|^{-1}\}.
\end{equation}

\section{The Sum $S_n$}

In order to investigate the sum $S_n$ we decompose the available ranges
for $r,r'$ and the corresponding zeros $\rho,\rho'$ into (overlapping)
ranges
\begin{equation}
\left\{\begin{array}{cc}R\le r\le RN^{\vp},&\;\;\;R'\le r'\le
R'N^{\vp},\\
T-1\le |\gamma|<TN^{\vp},&\;\;\; T'-1\le |\gamma'|<T'N^{\vp}.
\end{array}\right.
\end{equation}
Clearly $O(1)$ such ranges suffice to cover all possibilities,
so it is enough to consider the contribution from a fixed range of the
above type.  Throughout this section we shall follow the convention that 
$\rho=1$ is to included amongst the `zeros' corresponding to the trivial 
character.

Let $N(\sigma,\chi,T)$ denote as usual, the number of zeros $\rho$ of 
$L(s,\chi)$, in the region $\beta\ge \sigma$, $|\gamma|\le T$, and let
$N(\sigma,r,T)$ be the sum of $N(\sigma,\chi,T)$ for all characters $\chi$
of conductor $r$.  Since
\[N^{\beta-1}=N^{3\vp-1/2}+\int_{1/2+3\vp}^{\beta}N^{\sigma-1}
(\log N)d\sigma\]
for $\beta\ge 1/2+3\vp$, we find that
\begin{equation}
\sum_{\rho}N^{\beta-1}\ll N^{6\vp-1/2}RT+I(r)\log N,
\end{equation}
where the sum is over zeros of $L(s,\chi)$ for all $\chi$ of conductor
$r$, subject to $T-1\le |\gamma|\le TN^{\vp}$, and were
\[I(r)=\int_{1/2+3\vp}^{1}N^{\sigma-1}N(\sigma,r,TN^{\vp})d\sigma.\]

In view of the minimum occuring in (5) it is convenient to set
\[m(R,T)=\min(1,\frac{P}{RT}).\]
We now insert (7) into (5) so that, for given $r,r'$, the range (6)
contributes to 
\[\sum_{\chi\!\!\!\pmod{r}}A_m(\chi)\]
a total
\begin{eqnarray}
&\ll& \phi(r)\frac{P}{m}N^{4\vp}+\frac{N^{1/2}}{\log N}m(R,T)N^{6\vp-1/2}RT
+N^{1/2}m(R,T)I(r)\nonumber\\
&\ll&   PN^{6\vp}+N^{1/2}m(R,T)I(r).
\end{eqnarray}
Similarly, for the double sum
\[\sum_{\chi\!\!\!\pmod{r}}\sum_{\chi'\!\!\!\pmod{r'}}A_m(\chi)A_m(\chi')\]
the contribution is
\begin{equation}
\begin{array}{ll}\ll & P^{2}N^{12\vp}+PN^{1/2+6\vp}m(R,T)I(r)\\
&{}+PN^{1/2+6\vp}m(R',T')I(r')+Nm(R,T)m(R',T')I(r)I(r').
\end{array}
\end{equation}
We then sum over $r,r'$ using the following lemma.

\begin{lemma}
Let
\[\max_{r\le R}N(\sigma,r,T)=N_1(R),\;\;\;
\max_{r'\le R'}N(\sigma',r',T')=N_1(R'),\]
and
\[\sum_{r\le R}N(\sigma,r,T)=N_2(R),\;\;\;
\sum_{r'\le R'}N(\sigma',r',T')=N_2(R').\]
In the notation of (4) we have
\begin{eqnarray}
\lefteqn{\sum_{r\le R}\sum_{r'\le R'}
N(\sigma,r,T)N(\sigma',r',T')(ss')^{\vp-1}}\hspace{2cm}\\
&\ll &
\{N_1(R)N_2(R)N_1(R')N_2(R')\}^{1/2+2\vp},\nonumber
\end{eqnarray}
for $1/2\le\sigma,\sigma'\le 1$.

Moreover, if
\[P\le N^{45/154-4\vp},\]
then
\begin{equation}
\sum_{r,r'}m(R,T)m(R',T')
N(\sigma,r,TN^{\vp})N(\sigma',r',T'N^{\vp})(ss')^{\vp-1},
\end{equation}
\begin{equation}\ll N^{(1-\vp)(1-\sigma)+(1-\vp)(1-\sigma')}
\end{equation}
for $1/2+3\vp\le\sigma,\sigma'\le 1$, where the summation is for 
$R\le r\le RN^{\vp}$ and $R'\le r'\le R'N^{\vp}$.
\end{lemma}

We shall prove this at the end of this section.  Henceforth we shall
assume that $P\le N^{45/154-4\vp}$.  

For suitable values 
of $\eta$ in the range
\begin{equation}
0\le\eta\le \log\log N
\end{equation}
we shall define $\cl{B}(\eta)$ to be the set of characters $\chi$ of
conductor $r\le P$, for which the function $L(s,\chi)$ has at least
one zero in the region
\[\beta>1-\frac{\eta}{\log N},\;\;\; |\gamma|\le N.\]
According to our earlier convention the trivial character is always in 
$\cl{B}(\eta)$.
Now, if we restrict attention to pairs $\chi,\chi'$ for which
$\chi\not\in\cl{B}(\eta)$ we have
\begin{eqnarray*}
\lefteqn{\sum_{R\le r\le RN^{\vp}}\sum_{R'\le r'\le R'N^{\vp}}
Nm(R,T)m(R',T')I(r)I(r')(ss')^{\vp-1}}\hspace{3cm}\\
&\ll & \int_{1/2+3\vp}^{1-\eta/\log N}\int_{1/2+3\vp}^{1}
N^{1-\vp(1-\sigma)-\vp(1-\sigma')}d\sigma' d\sigma\\
&\ll & N^{1-\vp\eta/\log N}(\log N)^{-2}\\
&=& e^{-\vp\eta}N(\log N)^{-2}.
\end{eqnarray*}
Terms for which $\chi\in\cl{B}(\eta)$ but 
$\chi'\not\in\cl{B}(\eta)$ may be handled similarly.
This concludes our discussion of the final term in (9) for the time being.

To handle the third term in (9) we use the zero density
estimate
\begin{equation}
\sum_{r\le R}N(\sigma,r,T)\ll  (R^2 T)^{\kappa(\sigma)(1-\sigma)},
\end{equation}
where
\begin{equation}
\kappa(\sigma)=\left\{\begin{array}{cc} \frac{3}{2-\sigma}+\vp, & 
\frac{1}{2}\le\sigma\le\frac{3}{4}\\ \frac{12}{5}+\vp, &
\frac{3}{4}\le\sigma\le 1.\end{array}\right.
\end{equation}
This follows from results of
Huxley [5], Jutila [7; Theorem 1] and Montgomery [14; Theorem 12.2].  
For each fixed
value of $r'$ we have
\begin{eqnarray*}
\sum_{r} (ss')^{\vp-1}&\le&
\sum_{h|r'} (r'/h)^{\vp-1}\sum_{s\le P/h}s^{\vp-1}\\
&\ll &\sum_{h|r'} (r'/h)^{\vp-1}(P/h)^{\vp}\\
&\ll & N^{\vp}.
\end{eqnarray*}
The contribution of the third term in (9) to $S_n$ is therefore
\[\ll  PN^{1/2+5\vp}m(R',T')\sum_{r'}I(r').\]
However the bound (14) shows that
\[m(R',T')\sum_{r'}N(\sigma,r',TN^{\vp})\ll 
\min\{1\,,\,\frac{P}{R'T'}\}
({R'}^2 N^{2\vp}T'N^{\vp})^{\kappa(\sigma)(1-\sigma)}.\]
Since
\[0\le \kappa(\sigma)(1-\sigma)\le 1\]
in the range $1/2+\vp\le\sigma\le 1$, this is
\[\ll (P^2 N^{3\vp})^{\kappa(\sigma)(1-\sigma)}.\]
Moreover, if $P\le N^{45/154-4\vp}$, then
\[(P^2 N^{3\vp})^{\kappa(\sigma)(1-\sigma)}N^{\sigma-1}\le N^{f(\sigma)}\]
with
\begin{eqnarray*}
f(\sigma)&=&(\frac{45}{77}\kappa(\sigma)-1)(1-\sigma)\\
&\le& (\frac{45}{77}\{\frac{12}{5}+\vp\}-1)(1-\sigma)\\
&\le& (\frac{31}{77}+\vp)(1-\sigma)\\
&\le& (\frac{31}{77}+\vp)\frac{1}{2}\\
&\le& \frac{31}{154}+\vp.
\end{eqnarray*}
It follows that the contribution of the third term in (9) to $S_n$ is
\[\ll  PN^{1/2+6\vp}.N^{31/154+\vp}\ll N^{1-\vp}.\]
The second term may of course be handled similarly.

Finally we deal with the first term of (9) which produces a contribution
to $S_n$ which is
\begin{eqnarray*}
&\ll & P^{2}N^{12\vp}\sum_{r,r'}(ss')^{\vp-1}\\
&\ll & P^{2}N^{12\vp}\sum_{ss'h\le P}(ss')^{\vp-1}\\
&\ll & P^{2}N^{12\vp}\sum_{ss'\le P}P(ss')^{\vp-2}\\
&\ll & P^{3}N^{12\vp}\\
&\ll & N^{1-\vp},
\end{eqnarray*}
for $P\le N^{45/154-4\vp}$.

We summarize our conclusions thus far as follows.
\begin{lemma}
If $P\le N^{45/154-4\vp}$ then
\[S_n\le
\sum_{\chi,\chi'\in\cl{B}(\eta)}A_m(\chi)A_m(\chi')m_n^{-1/3}
+O(e^{-\vp\eta}N(\log N)^{-2}).\]
\end{lemma}

To handle the characters in $\cl{B}(\eta)$ we use the zero-density
estimate
\begin{equation}
N(\sigma,r,T)\ll  (rT)^{\kappa(\sigma)(1-\sigma)},
\end{equation}
with $\kappa(\sigma)$ given by (15).
This also follows from work of 
Huxley [5], Jutila [7; Theorem 1] and Montgomery [14; Theorem 12.1].  Thus
\begin{eqnarray*}
m(R,T)N(\sigma,r,TN^{\vp})&\ll &
\max\{1\,,\,\frac{P}{RT}\}(rTN^{\vp})^{\kappa(\sigma)(1-\sigma)}\\
&\ll &(PN^{2\vp})^{\kappa(\sigma)(1-\sigma)}\\
&\ll &(PN^{2\vp})^{(12/5+\vp)(1-\sigma)}\\
&\ll &N^{(1-\vp)(1-\sigma)}
\end{eqnarray*}
for $P\le N^{45/154-4\vp}$.  We deduce that
\[m(R,T)I(r)\ll (\log N)^{-1}.\]
It follows from (8) that
\[A_m(\chi)\ll N^{1/2}(\log N)^{-1}.\]
We also note that
\[\#\cl{B}(\eta)\ll\sum_{r}N(1-\frac{\eta}{\log N},r,N)\ll
(P^2 N)^{3\eta/\log N}\ll e^{6\eta},\]
by (14), since $\kappa(\sigma)\le 3$ for all $\sigma$.  
We therefore have the following facts.
\begin{lemma}
If $\chi\in\cl{B}(\eta)$, we have $A_m(\chi)\ll N^{1/2}(\log
N)^{-1}$.  Moreover, we have $\#\cl{B}(\eta)\ll e^{6\eta}$.
\end{lemma}

We end this section by establishing Lemma 3.  
We shall suppose, as we may by the symmetry, that
\begin{equation}
N_2(R)N_1(R')\le N_2(R')N_1(R).
\end{equation}
Let $U\ge 1$ be a parameter whose value will be assigned in due course,
see (18).  For those terms of the sum (10) in which
$ss'\ge U$ we plainly have a total
\[\le \sum_{r\le R}\sum_{r'\le R'}N(\sigma,r,T)N(\sigma',r',T')U^{\vp-1}
\ll N_2(R)N_2(R')U^{\vp-1}.\]
On the other hand, when $ss'<U$
we observe that, for fixed $s,s'$ we have
\begin{eqnarray*}
\sum_{h}N(\sigma,hs,T)N(\sigma',hs',T')&\ll&
\sum_{h}N(\sigma,hs,T)N_1(R')\\
&\ll& \sum_{r}N(\sigma,r,T)N_1(R')\\
&\ll& N_2(R)N_1(R').\\
\end{eqnarray*}
On summing over $s$ and $s'$ we therefore obtain a total
\[\ll N_2(R)N_1(R')\sum_{ss'\le U}(ss')^{\vp-1}\ll 
N_2(R)N_1(R')U^{2\vp}.\]
It follows that the sum (10) is
\[\ll N_2(R)\{N_2(R')U^{2\vp-1}+N_1(R')U^{2\vp}\}.\]
We therefore choose
\begin{equation}
U=N_2(R')/N_1(R'),
\end{equation}
whence the sum (10) is
\begin{eqnarray*}
&\ll& N_2(R)N_1(R')U^{2\vp}\\
&\ll& N_2(R)N_1(R')\{N_1(R)N_2(R)N_1(R')N_2(R')\}^{2\vp}\\
&\ll& \{N_2(R)N_1(R')N_2(R')N_1(R)\}^{1/2}
\{N_1(R)N_2(R)N_1(R')N_2(R')\}^{2\vp}\\
\end{eqnarray*}
in view of (17).  This produces the required bound.

To establish (12) we shall bound $N_1(R)$ and $N_1(R')$ using (16).
Moreover to handle $N_2(R)$ and $N_2(R')$ we shall use the estimate
\[\sum_{r\le R}N(\sigma,r,T)\ll  \left\{\begin{array}{cc}
(R^2 T)^{\kappa(\sigma)(1-\sigma)},\; &\; 
\frac{1}{2}+\vp\le\sigma\le \frac{23}{38}\\
(R^2 T^{6/5})^{\lambda(1-\sigma)},\; &\;
\frac{23}{38}<\sigma\le 1,\end{array}\right.\]
where
\[\lambda=\frac{20}{9}+\vp.\]
This follows from (14) and (15) along with Heath-Brown [4; Theorem 2] and
Jutila [7; Theorem 1].

We now see that the sum (11) may be estimated as
\begin{equation}
\ll m(R,T)R^aT^c.m(R',T'){R'}^b{T'}^d. N^{e},
\end{equation}
say, where
\[a=\left\{\begin{array}{cc} 3\kappa(\sigma)(1-\sigma)(\frac{1}{2}+2\vp),
\; &\; 
\frac{1}{2}+3\vp\le\sigma\le \frac{23}{38}\\
\{\kappa(\sigma)+2\lambda\}(1-\sigma)(\frac{1}{2}+2\vp),\; &\;
\frac{23}{38}<\sigma\le 1,\end{array}\right.\]
and 
\[c=\left\{\begin{array}{cc} 2\kappa(\sigma)(1-\sigma)(\frac{1}{2}+2\vp)
,\; &\; 
\frac{1}{2}+3\vp\le\sigma\le \frac{23}{38}\\
\{\kappa(\sigma)+6\lambda/5\}(1-\sigma)(\frac{1}{2}+2\vp),\; &\;
\frac{23}{38}<\sigma\le 1,\end{array}\right.\]
and similarly for $b$ and $d$.  Moreover we may take
\[e=6\vp(1-\sigma)+6\vp(1-\sigma').\]
It therefore follows that $0\le c,d< 1$,
whence (19) is maximal for $T=P/R$ and $T'=P/R'$.  Similarly we have
$a\ge c$ and $b\ge d$.  Thus, after substituting $T=P/R$ and
$T'=P/R'$ in (19), the resulting expression is increasing with respect
to $R$ and $R'$, and hence is maximal when $R=R'=P$.  We therefore see
that (20) is 
\[\ll P^{a+b}N^e.\]
Finally one can check that
\[(\frac{45}{154}-4\vp)a\le (1-7\vp)(1-\sigma),\]
and similarly for $b$.  This suffices to establish the
bound (12) for $P\le N^{45/154-4\vp}$.

\section{Summation Over Powers of 2}

In this section we consider the major arc integral
\[\int_{{\mathfrak M}}S(\alpha)^2 T(\alpha)^K e(-\alpha N)d\alpha,\]
where we now assume $N$ to be even.
According to Lemmas 2 and 4 we have
\begin{eqnarray}
\int_{{\mathfrak M}}S(\alpha)^2 T(\alpha)^K e(-\alpha N)d\alpha&=&
\Sigma_0+O(e^{-\vp\eta}N(\log N)^{-2}\Sigma_1)\nonumber\\
&&\hspace{1cm}+O(N(\log N)^{-2}\Sigma_2),
\end{eqnarray}
where 
\[\Sigma_0=\sum_{n} J_n(1,1),\]
\[\Sigma_1=\sum_{n}\frac{n}{\phi(n)}\]
and
\[\Sigma_2=\sum_{\chi,\chi'\in\cl{B}(\eta)}\sum_{n}
\frac{n}{\phi(n)}m_n^{-1/3}.\]
In each case the sum over $n$ is for values
\begin{equation}
n=N-\sum_{j=1}^{K}2^{\nu_j}.
\end{equation}

We begin by considering the main term $\Sigma_0$.  We put
\[T(\beta)=\sum_{\vp N<m\le N}\frac{e(\beta m)}{\log m}\]
and
\[R(\beta)=S(\beta)-T(\beta).\]
We also set
\[||R||=\int_{-P/N}^{P/N}|R(\beta)|^2 d\beta\]
and
\[J(n)=\twosum{\vp<m_1,m_2<N}{m_1+m_2=n}(\log m_1)^{-1}(\log m_2)^{-1}.\]
Then, as in Gallagher [3; (11)], we have
\begin{eqnarray}
J_n(1,1)&=&J(n)\cl{S}(n)
+O(N(\log N)^{-2}\frac{n}{\phi(n)}d(n)\frac{\log P}{P})\nonumber\\
&&\hspace{1cm}+O(\frac{n}{\phi(n)}\{N^{1/2}(\log N)^{-1}||R||+||R||^2\}),
\end{eqnarray}
where
\[\cl{S}(n)=\prod_{p|n}(\frac{p}{p-1})\prod_{p\nd n}(1-\frac{1}{(p-1)^2}).\]
In analogy to (5) we have
\[||R||\ll PN^{4\vp}+\frac{N^{1/2}}{\log N}
\sum_{\rho}N^{\beta-1}\min\{1\,,\,P|\gamma|^{-1}\},\]
where the sum over $\rho$ is for zeros of $\zeta(s)$ in the
region
\[\beta\ge\frac{1}{2}+3\vp,\;\;\;|\gamma|\le N.\]
We split the range for $|\gamma|$ into $O(1)$ overlapping
intervals
\[T-1\le |\gamma|\le TN^{\vp},\]
and find, as in (8) that each range contributes
\[\ll PN^{4\vp}+N^{1/2}\min\{1\,,\,\frac{P}{T}\}
\{N^{6\vp-1/2}T+\int_{1/2+3\vp}^{1}N^{\sigma-1}N(\sigma,1,TN^{\vp})d\sigma\}\]
to $||R||$.  Using the case $R=1$ of (14), together with Vinogradov's
zero-free region
\[\sigma\ge 1-\frac{c_0}{(\log T)^{3/4}(\log\log T)^{3/4}}\]
(see Titchmarsh [16; (6.15.1)]), we find that this gives
\[||R||\ll N^{1/2}(\log N)^{-10},\]
say, for $P\le N^{45/154-4\vp}$.  The error terms in (22) are
therefore $O(N(\log N)^{-9})$.

We also note that
\begin{eqnarray*}
J(n)&=&(\log N)^{-2}\#\{m_1,m_2:\vp N<m_1,m_2\le N,\,m_1+m_2=n\}\\
&&\hspace{3cm}+O(N(\log N)^{-3})\\
&=&(\log N)^{-2}R(n)+O (N(\log N)^{-3}),
\end{eqnarray*}
where
\[R(n)=\left\{\begin{array}{cc} 2N-n,& (1+\vp)N\le n\le 2N,\\
n-2\vp N,& 2\vp N\le n\le (1+\vp)N,\\
0, & \mbox{otherwise}.\end{array}\right.\]
In particular, we have $R(N-m)=(1-2\vp)N(\log N)^{-2}+O(m(\log
N)^{-2})$ for $1\le m\le N$.  Since 
\[\cl{S}(n)\ll\frac{n}{\phi(n)}\ll\log\log N, \]
we find, on taking $n$ of the form (21), that
\[\sum_{n}J(n)\cl{S}(n)=(1-2\vp)N(\log
N)^{-2}\sum_{n}\cl{S}(n)+O(N(\log N)^{K-5/2})\]
for $K\ge 2$, whence
\[\Sigma_0=(1-2\vp)N(\log N)^{-2}\sum_{n}\cl{S}(n)+
O(N(\log N)^{K-5/2}).\]

Since the numbers $n$ are all even, we have
\[\cl{S}(n)=2C_0\prod_{p|n,
p\not=2}\frac{p-1}{p-2}=2C_0\sum_{d|n}k(d),\]
where
\begin{equation}
C_0=\prod_{p\not=2}(1-\frac{1}{(p-1)^2})
\end{equation}
and $k(d)$ is the multiplicative function defined by taking 
\begin{equation}
k(p^e)=\left\{\begin{array}{cc}
0,\; & p=2\;\mbox{or}\;e\ge 2,\\
(p-2)^{-1},\;& \mbox{otherwise.}
\end{array}\right.
\end{equation}
For any odd integer $d$ we shall define $\ep(d)$ to be the order of 2
in the multiplicative group modulo $d$, and we shall set
\[H(d;N,K)=\#\{(\nu_1,\ldots,\nu_K): 1\le\nu_i\le\ep(d),\,
d|N-\sum 2^{\nu_i}\}.\]
Then for any fixed $D$ we have
\begin{eqnarray*}
\sum_{n}\cl{S}(n)&=&2C_0\sum_{d}k(d)\#\{n:d|n\}\\
&\ge&2C_0\sum_{d\le D}k(d)\#\{n:d|n\}\\
&\ge&2C_0\sum_{d\le D}k(d)H(d;N,K)[L/\ep(d)]^K\\
&\ge&\{1+O((\log N)^{-1})\}2C_0L^K 
\sum_{d\le D}k(d)H(d;N,K)\ep(d)^{-K}.
\end{eqnarray*}
We shall take $D=5$.
We trivially have $\ep(1)=1$ and $H(1;N,K)=1$ for all $N$ and $K$.
When $d=3$ or $d=5$ the powers of 2 run over all non-zero residues
modulo $d$, and it is an easy exercise to check that
\[H(d;N,K)=\left\{\begin{array}{cc} \frac{1}{d}\{(d-1)^K-(-1)^K\},
& d\nd N\\ \frac{1}{d}\{(d-1)^K+(-1)^K (d-1)\}, & d|N.\end{array}\right.\]
Thus if $K\ge 7$ we have
\[H(3;N,K)\ep(3)^{-K}\ge \frac{1}{3}(1-2^{-6})\]
and
\[H(5;N,K)\ep(5)^{-K}\ge \frac{1}{5}(1-4^{-6}),\]
whence
\[2\sum_{d\le D}k(d)H(d;N,K)\ep(d)^{-K}\ge 2.7895\]
for any choice of $N$.  We therefore conclude that
\begin{equation}
\Sigma_0\ge 2.7895(1-2\varpi)C_0N(\log N)^{-2}L^K+O(N(\log N)^{K-5/2}),
\end{equation}
providing that $K\ge 9$.

To bound $\Sigma_1$ we note that
\[\frac{n}{\phi(n)}\ll\prod_{p|n,\,p\not=2}(1+\frac{1}{p})=
\sum_{q|n,\,2\nd q}\frac{\mu^2(q)}{q}.\]
We deduce that
\[\Sigma_1\ll \sum_{q\le N,\,2\nd q}\frac{\mu^2(q)}{q}\#\{n:\,q|n\}.\]
However, if $q$ is odd, then
\[\#\{\nu:0\le\nu\le L,\,2^{\nu}\equiv m\!\!\!\pmod{q}\}\ll
1+\frac{L}{\ep(q)}.\]
It follows that
\[\#\{n:\,q|n\}\ll L^{K-1}+L^K /\ep(q),\]
whence
\[\Sigma_{1}\ll (\log N)^K+
(\log N)^K\sum_{q\le N,\,2\nd q}\frac{\mu^2(q)}{q\ep(q)}.\]
To bound the final sum we call on the following simple result of Gallagher
[3; Lemma 4]
\begin{lemma}
We have
\[\sum_{\ep(q)\le x}\frac{\mu^2(q)}{\phi^2(q)}q\ll\log x.\]
\end{lemma}
From this we deduce that
\begin{equation}
\sum_{x/2<\ep(q)\le x}\frac{\mu^2(q)}{q\ep(q)}\ll\frac{\log x}{x}.
\end{equation}
We take $x$ to run over powers of $2$ and sum the resulting bounds to
deduce that
\[\sum_{q\le N,\,2\nd q}\frac{\mu^2(q)}{q\ep(q)}\ll 1,\]
and hence that 
\begin{equation}
\Sigma_{1}\ll (\log N)^K.
\end{equation}

Turning now to $\Sigma_2$, we fix a particular pair of characters
$\chi,\chi'\in\cl{B}(\eta)$, and investigate
\[\sum_{n}\frac{n}{\phi(n)}m_n^{-1/3}=\Sigma_2(\chi,\chi'),\]
say.  Let $m=[r,r']$ as usual, and write $m=2^{\mu}f$, with $f$ odd. 
Put $g=(f,n)$ so that
\begin{equation}
m_n\ge f_n=f/g,
\end{equation}
and consider
\[\sum_{g|n}\frac{n}{\phi(n)}.\]
As before we have
\[\frac{n}{\phi(n)}\ll\sum_{q|n,\,2\nd q}\frac{\mu^2(q)}{q}.\]
Terms $q$ with $q\ge d(n)$ can contribute at most $1$ in total, so that
in fact
\[\frac{n}{\phi(n)}\ll\sum_{q|n,\,2\nd q, q\le d(n)}\frac{\mu^2(q)}{q}.\]
Thus, if
\[D=\max_{1\le n\le N}d(n),\]
we deduce as before that
\begin{eqnarray*}
\sum_{g|n}\frac{n}{\phi(n)}&\ll&\sum_{q\le D,\,2\nd q}
\frac{\mu^2(q)}{q}\#\{n:\,[g,q]|n\}\\
&\ll &\sum_{q\le D,\,2\nd q}\frac{\mu^2(q)}{q}\{(\log
N)^{K-1}+\frac{(\log N)^K}{\ep([g,q])}\}.
\end{eqnarray*}
Here we note that
\[\sum_{q\le D}q^{-1}\ll \log D\ll\frac{\log N}{\log\log N}.\]
To deal with the remaining terms let $\xi$ be a positive parameter.  Then
\begin{eqnarray*}
\sum_{\ep(q)>\xi}\frac{\mu^2(q)}{q\ep([g,q])}&\le &
\sum_{\ep(q)>\xi}\frac{\mu^2(q)}{q\ep(q)}\\
&\ll& \frac{\log{\xi}}{\xi},
\end{eqnarray*}
by (26).  If $\ep(q)\le \xi$ we note that 
\begin{equation}
q\le 2^{\ep(q)}-1,\;\; \mbox{for}\;\;q>1,
\end{equation}
so that $q\le 2^{\xi}$.  Thus
\begin{eqnarray*}
\sum_{\ep(q)\le\xi}\frac{\mu^2(q)}{q\ep([g,q])}&\le &
\sum_{q\le 2^{\xi}}\frac{\mu^2(q)}{q\ep(g)}\\
&\le & \frac{\xi}{\ep(g)}.
\end{eqnarray*}
On choosing $\xi=\sqrt{\ep(g)}$ we therefore conclude that
\[\sum_{2\nd q}\frac{\mu^2(q)}{q\ep([g,q])}\ll
\frac{\log\ep(g)}{\sqrt{\ep(g)}},\]
and hence that
\[\sum_{g|n}\frac{n}{\phi(n)}\ll 
(\log N)^K \{(\log\log N)^{-1}+\ep(g)^{-1/3}\}.\]
It follows from (29) that $\ep(g)\gg\log g$, and we now conclude that
\[\sum_{g|n}\frac{n}{\phi(n)}\ll 
(\log N)^K \{(\log\log N)^{-1}+(\log g)^{-1/3}\}.\]

We now observe from (28) that
\[\Sigma_2(\chi,\chi')\le
\sum_{n}\frac{n}{\phi(n)}(\frac{f}{(f,n)})^{-1/3}.\]
Let $\tau\ge 1$ be a parameter to be fixed in due course. Then terms 
in which $(f,n)\le f/\tau$ contribute 
\[\le \tau^{-1/3}\sum_{n}\frac{n}{\phi(n)}=\tau^{-1/3}\Sigma_1
\ll \tau^{-1/3}(\log N)^K,\]
by (27).  The remaining terms contribute
\begin{eqnarray*}
&\le&\sum_{g|f,\,g\ge f/\tau}(f/g)^{-1/3}\sum_{g|n}\frac{n}{\phi(n)}\\
&\ll& \sum_{g|f,\,g\ge f/\tau}(f/g)^{-1/3}
(\log N)^K \{(\log\log N)^{-1}+(\log g)^{-1/3}\}\\
&\ll& \sum_{g|f,\,g\ge f/\tau}
(\log N)^K \{(\log\log N)^{-1}+(\log f)^{-1/3}\}\\
&\ll& \sum_{j|f,\,j\le \tau}
(\log N)^K \{(\log\log N)^{-1}+(\log f)^{-1/3}\}\\
&\ll& \tau(\log N)^K \{(\log\log N)^{-1}+(\log f)^{-1/3}\}.
\end{eqnarray*}
We deduce that
\[\Sigma_2(\chi,\chi')\ll\tau^{-1/3}(\log N)^K+
\tau(\log N)^K \{(\log\log N)^{-1}+(\log f)^{-1/3}\}.\]
We therefore choose
\[\tau=\{(\log\log N)^{-1}+(\log f)^{-1/3}\}^{-3/4},\]
whence
\begin{equation}
 \Sigma_2(\chi,\chi')\ll (\log N)^K\{(\log\log N)^{-1/4}+(\log f)^{-1/12}\}.
\end{equation}

In order to bound $f$ from below we note that, since $\chi,\chi'$ are 
not both trivial, we may suppose 
that $\chi$, say, is non-trivial.  We then
use a result of Iwaniec [6;~Theorem~2].  
This shows that if $L(\beta+i\gamma,\chi)=0$, with
$|\gamma|\le N$, and $\chi$ of conductor $r\le N$, then either $\chi$ is
real, or
\[1-\beta\gg \{\log d+(\log N\log\log N)^{3/4}\}^{-1},\]
where $d$ is the product of the distinct prime factors of $r$.  In our
application we clearly have $f\ge d/2$, so that if $\chi$, say, is in
$\cl{B}(\eta)$ we must have
\[\frac{\eta}{\log N}\gg \{\log f+(\log N\log\log N)^{3/4}\}^{-1}\]
if $\chi$ is not real.
Thus, if we insist that $\eta\le (\log N)^{1/5}$ it follows that either
\[\log f\gg\eta^{-1}\log N\gg (\log N)^{4/5},\]
or $\chi$ is real.  Of course if $\chi$ is real we will have $16\nd
r$, whence $f\gg r$.  Moreover we will also have 
\[(\log N)^{4/5}\gg\frac{\eta}{\log N}\gg 1-\beta\gg r^{\vp-1/2},\]
so that $f\gg r\gg(\log N)^{3/2}$.  Thus in either case we find that
$\log f\gg\log\log N$, so that (30) yields
\[ \Sigma_2(\chi,\chi')\ll (\log N)^K(\log\log N)^{-1/12}.\]
In view of the bound for $\#\cl{B}(\eta)$ given in Lemma 5, we
conclude that
\begin{equation}
\Sigma_2\ll e^{12\eta}(\log N)^K(\log\log N)^{-1/12}.
\end{equation}

We may now insert the bounds (25), (27) and (31) into (20) to deduce that
\begin{eqnarray*}
\int_{{\mathfrak M}}S(\alpha)^2 T(\alpha)^K e(-\alpha N)d\alpha
&\ge &2.7895(1-2\varpi)C_0N(\log N)^{-2}L^K\\
&&\hspace{2mm}{}+O(N(\log N)^{K-5/2})\\
&&\hspace{4mm}{}+O(e^{-\vp\eta}N(\log N)^{K-2})\\
&&\hspace{6mm}{}+O(e^{12\eta}N(\log N)^{K-2}(\log\log N)^{-1/12}).
\end{eqnarray*}
We therefore define $\eta$ by taking
\[e^{\eta}=(\log\log N)^{1/145},\]
so that $\eta$ satisfies the condition (13), and conclude as follows.
\begin{lemma}
If $p\le N^{45/154-4\vp}$ and $K\ge 9$ we have
\[\int_{{\mathfrak M}}S(\alpha)^2 T(\alpha)^K e(-\alpha N)d\alpha
\ge 2.7895(1-3\varpi)C_0N(\log 2)^{-2}L^{K-2}\]
for large enough $N$.
\end{lemma}

\section{A Mean Square Estimate}

In this section we shall estimate the mean square
\[J({\mathfrak m})=\int_{{\mathfrak m}}|S(\alpha)T(\alpha)|^2 d\alpha,\]
where ${\mathfrak m}=[0,1]\setminus{\mathfrak M}$ is the set of minor arcs.
Instead of this integral, previous researchers have worked with the
larger integral
\[J=\int_{0}^1 |S(\alpha)T(\alpha)|^2 d\alpha.\]
Thus it was shown by Li [9; Lemma 6], building on work of Liu, Liu and 
Wang [13; Lemma 4] that
\[J\le (24.95+o(1))\frac{C_0}{\log^2 2}N,\]
In this section we shall improve on this bound, and give a lower bound for
the corresponding major arc integral
\[J({\mathfrak M})=\int_{{\mathfrak M}}|S(\alpha)T(\alpha)|^2 d\alpha.\]
By subtraction we shall then obtain our bound for $J({\mathfrak m})$.

We begin by observing that
\[J=\sum_{\mu,\nu\le L}r(2^{\mu}-2^{\nu}),\]
where
\[r(n)=\#\{\varpi N<p_1,p_2\le N: n=p_1-p_2\}.\]
Moreover, by  Theorem 3 of Chen [2] we have
\[r(n)\le C_0 C_1 h(n)\frac{N}{(\log N)^2},\]
for $n\not=0$ and $N$ sufficiently large, where $C_0$ is given by (23),
\begin{equation}
C_1=7.8342,
\end{equation}
and
\[h(n)=\prod_{p|n,\,p>2}(\frac{p-1}{p-2}).\]
Observe that our notation for the constants that occur differs from that
used by Liu, Liu and Wang, and by Li.  Since
$h(2^{\mu}-2^{\nu})=h(2^{\mu-\nu}-1)$ for $\mu>\nu$ we conclude, as in 
Liu, Liu and Wang [13; \S 3] and Li [9; \S 4] that
\begin{equation}
\sum_{\mu\not=\nu\le L}r(2^{\mu}-2^{\nu})\le 2C_0 C_1\frac{N}{(\log N)^2}
\sum_{1\le l\le L}(L-l)h(2^l-1),
\end{equation}
while the contribution for $\mu=\nu$ is $L\pi(N)-L\pi(\varpi N)\le
LN(\log N)^{-1}$, for large $N$.  Now
\[h(n)=\sum_{d|n}k(d),\]
where $k(d)$ is the multiplicative function defined in (24).  Thus
\begin{eqnarray*}
\sum_{1\le j\le J}h(2^j-1)&=&
\sum_{d=1}^\infty k(d)\#\{j\le J: d|2^j-1\}\\
&=&\sum_{d=1}^\infty k(d)[\frac{J}{\ep(d)}].
\end{eqnarray*}
However $[\theta]=\theta+O(\theta^{1/2})$ for any real $\theta>0$,
whence
\begin{equation}
\sum_{1\le j\le J}h(2^j-1)=C_2 J+O(J^{1/2})
\end{equation}
with
\begin{equation}
C_2=\sum_{d=1}^\infty \frac{k(d)}{\ep(d)}.
\end{equation}
Here we use the observation that the sum
\[\sum_{d=1}^\infty \frac{k(d)}{\ep(d)^{1/2}}\]
is convergent, since Lemma 6 implies that
\begin{equation}
\sum_{x/2<\ep(d)\le x}\frac{k(d)}{\ep(d)^{1/2}}\ll
x^{-1/2}\sum_{x/2<\ep(d)\le x}\frac{\mu^2(d)d}{\phi^2(d)}\ll
\frac{\log x}{x^{1/2}}
\end{equation}
for any $x\geq 2$.

We may now use partial summation in conjunction with (34) to deduce that
\[\sum_{1\le l\le L}(L-l)h(2^l-1)=C_2\frac{L^2}{2}+O(L^{3/2}),\]
Thus, using (33) we reach the following result.
\begin{lemma}
We have
\[J\le \{\frac{C_0 C_1 C_2}{\log^2 2}+\frac{1}{\log 2}+o(1)\}N,\]
with the constants given by (23), (32) and (35).
\end{lemma}

We now turn to the integral $J({\mathfrak M})$.
According to Lemma 3.1 of Vaughan [17], if 
\[|\alpha-\frac{a}{q}|\le\frac{\log x}{x},\;\;\;(a,q)=1,\]
and $q\le 2\log x$, we have
\[\sum_{p\le x}e(\alpha p)\log p=\frac{\mu(q)}{\phi(q)}v(\alpha-\frac{a}{q})
+O(x(\log x)^{-3}),\]
with
\[v(\beta)=\sum_{m\le x}e(\beta m).\]
It follows by partial summation that
\[S(\alpha)=\frac{\mu(q)}{\phi(q)}w(\alpha-\frac{a}{q})
+O(N(\log N)^{-4}),\]
with
\[w(\beta)=\sum_{\vp N<m\le N}\frac{e(\beta m)}{\log m},\]
providing that 
\begin{equation}
|\alpha-\frac{a}{q}|\le\frac{\log N}{N},\;\;\;(a,q)=1
\end{equation}
and $q\le\log N$.  Then if $\mathfrak{a}$ denotes the set of $\alpha\in[0,1]$
for which such $a,q$ exist, we easily compute that
\begin{eqnarray*}
 J(\mathfrak{M})&\ge&J(\mathfrak{a})\\
&=& \int_{{\mathfrak a}}|\frac{\mu(q)}{\phi(q)}w(\alpha-\frac{a}{q})
T(\alpha)|^2 d\alpha+O(N(\log N)^{-1}),
\end{eqnarray*}
where, for each $\alpha\in\mathfrak{a}$, we have taken $a/q$ to be the
unique rational satisfying (37).
By partial summation we have
\[w(\beta)\ll (||\beta||\log N)^{-1},\]
whence
\[\int_{-(\log N)/N}^{(\log N)/N}|w(\beta)T(\frac{a}{q}+\beta)|^2 d\beta
=\int_{-1/2}^{1/2}|w(\beta)T(\frac{a}{q}+\beta)|^2 d\beta+O(N(\log N)^{-1}).
\]
It follows that
\[J(\mathfrak{a})=
\sum_{q\le\log N}\sum_{(a,q)=1}\frac{\mu^2(q)}{\phi^2(q)}
\int_{0}^{1}|w(\beta)T(\frac{a}{q}+\beta)|^2 d\beta+
O(N(\log N)^{-1}\log\log N).\]
The integral on the right is
\[\sum_{0\le\mu,\nu\le L}e(a(2^{\mu}-2^{\nu})/q)S(2^{\mu}-2^{\nu}),\]
where
\begin{eqnarray*}
S(n)&=&
\twosum{\vp N<m_1,m_2\le N}{m_1-m_2=n}(\log m_1)^{-1}(\log m_2)^{-1}\\
&=&(\log N)^{-2}\#\{m_1,m_2:\vp N<m_1,m_2\le N,\,m_1-m_2=n\}\\
&&\hspace{3cm}+O(N(\log N)^{-3})\\
&=&(\log N)^{-2}\max\{N(1-\vp)-|n|\,,\,0\}+O (N(\log N)^{-3}).
\end{eqnarray*}
Thus
\begin{equation}
S(n)=(1-\vp)N(\log N)^{-2}+O(|n|(\log N)^{-2})+O (N(\log N)^{-3})
\end{equation}
for $n\ll N$.
On summing over $a$ we now obtain
\[J(\mathfrak{a})=
\sum_{0\le\mu,\nu\le L}\sum_{q\le\log N}\frac{\mu^2(q)}{\phi^2(q)}
c_q(2^{\mu}-2^{\nu})S(2^{\mu}-2^{\nu})+O(N(\log N)^{-1}\log\log N),\]
where $c_q(n)$ is the Ramanujan sum.  When $q$ is square-free we have
$c_q(n)=\mu(q)\mu((q,n))\phi((q,n))$.  Thus the error terms in (38)
make a total contribution $O(N(\log N)^{-1}\log\log N)$ to 
$J(\mathfrak{a})$.  Moreover
\[\mu^2(q)c_q(n)=\mu(q)\sum_{d|(q,n)}\mu(d)d,\]
whence
\[\sum_{0\le\mu,\nu\le L}\mu^2(q)c_q(n)=\mu(q)\sum_{d|q}\mu(d)d
\#\{\mu,\nu:\,1\le\mu,\nu\le L,\,d|2^{\mu}-2^{\nu}\}.\]
If $d$ is odd we have
\[\#\{\mu,\nu:\,1\le\mu,\nu\le L,\,d|2^{\mu}-2^{\nu}\}=L^2\ep(d)^{-1}+O(L),\]
while if $d$ is even, of the form $2e$ with $e$ odd, we have
\[\#\{\mu,\nu:\,1\le\mu,\nu\le L,\,d|2^{\mu}-2^{\nu}\}=L^2\ep(e)^{-1}+O(L).\]
The error terms contribute $O(N(\log N)^{-1}\log\log N)$ to 
$J(\mathfrak{a})$, by (38), so that
\[J(\mathfrak{a})=\frac{(1-\vp)N}{(\log N)^2}L^2
\sum_{q\le\log N}\frac{\mu(q)}{\phi^2(q)}\sum_{d|q}\mu(d)d\ep(d)^{-1}
+O(N(\log N)^{-1}\log\log N),\]
where $\ep(d)$ is to be interpreted as $\ep(e)$ when $d=2e$.
Now
\begin{eqnarray}
\sum_{q\le\log N}\frac{\mu(q)}{\phi^2(q)}\sum_{d|q}\frac{\mu(d)d}{\ep(d)}&=&
\sum_{d\le\log N}\frac{\mu(d)d}{\ep(d)}
\twosum{q\le\log N}{d|q}\frac{\mu(q)}{\phi^2(q)}\nonumber\\
&=&\sum_{d\le\log N}\frac{\mu(d)d}{\ep(d)}
\sum_{j\le (\log N)/d}\frac{\mu(jd)}{\phi^2(jd)}\nonumber\\
&=&\sum_{d\le\log N}\frac{\mu^2(d)d}{\ep(d)\phi^2(d)}
\twosum{j\le (\log N)/d}{(j,d)=1}\frac{\mu(j)}{\phi^2(j)}\nonumber\\
&=&\sum_{d\le\log N}\frac{\mu^2(d)d}{\ep(d)\phi^2(d)}
\{\twosum{j=1}{(j,d)=1}^{\infty}\frac{\mu(j)}{\phi^2(j)}
+O(\frac{d}{\log N})\}\nonumber\\
&=&\sum_{d\le\log N}\frac{\mu^2(d)d}{\ep(d)\phi^2(d)}
\prod_{p\nd d}\{1-(p-1)^{-2}\}\nonumber\\
&&\hspace{1cm}+
O((\log N)^{-1}\sum_{d\le\log N}\frac{\mu^2(d)d^2}{\ep(d)\phi^2(d)}).
\end{eqnarray}
If $d=2e$ with $e$ odd, we have
\[\frac{\mu^2(d)d}{\ep(d)\phi^2(d)}
\prod_{p\nd d}\{1-(p-1)^{-2}\}=2C_0 k(e)/\ep(d),\]
while if $d$ is odd we have
\[\prod_{p\nd d}\{1-(p-1)^{-2}\}=0,\]
since the factor with $p=2$ vanishes.  Moreover
\[\sum_{d\gg\log N}\frac{k(d)}{\ep(d)}\ll
\frac{\log N}{\log\log N}\]
by Lemma 6, applied as in (36).  
The leading term in (39) is therefore $2C_0 C_2+o(1)$, with 
$C_0$ and $C_2$ as in (23) and (35).

To bound the error term we use Lemma 6, which shows that
\[\twosum{X<d\le 2X}{x<\ep(d)\le 2x}\frac{\mu^2(d)d^2}{\ep(d)\phi^2(d)}
\ll\frac{X\log x}{x}.\]
According to (29) we must have $x\gg\log X$, so on summing as $x$ runs
over powers of $2$ we obtain
\[\sum_{X<d\le 2X}\frac{\mu^2(d)d^2}{\ep(d)\phi^2(d)}
\ll\frac{X\log\log X}{\log X}.\]
Now, summing as $X$ runs over powers of $2$ we conclude that
\[\sum_{d\le\log N}\frac{\mu^2(d)d^2}{\ep(d)\phi^2(d)}\ll 
\frac{(\log N)(\log\log\log N)}{\log\log N}.\]

We may therefore summarize our results as follows.
\begin{lemma}
We have
\[ J(\mathfrak{M})\ge \{\frac{2(1-\vp)C_0 C_2}{\log^2 2}+o(1)\}N,\]
and hence
\[J(\mathfrak{m})\le \{\frac{C_0(C_1-2+2\vp)C_2}{\log^2 2}
+\frac{1}{\log 2}+o(1)\}N,\]
by Lemma 8.
\end{lemma}

It remains to compute the constants.  We readily find
\[\prod_{2<p\le 200000}(1-(p-1)^{-2})=0.6601...\]
Since
\[\prod_{p>K}(1-(p-1)^{-2})\ge\prod_{n=K}^{\infty}(1-n^{-2})
=1-K^{-1},\]
we deduce that
\begin{equation}
C_0\ge 0.999995\times0.6601\ge 0.66.
\end{equation}
However the estimation of $C_2$ is more difficult.   We set
\[m=\prod_{e\le x}(2^e-1)\]
and
\[s(x)=\sum_{\ep(d)\le x}k(d),\]
whence
\begin{eqnarray*}
s(x)&\le&\sum_{d|m}k(d)\\
&=& h(m)\\
&=&\prod_{p|m,\,p>2}(\frac{p-1}{p-2})\\
&\le & \prod_{p>2}(\frac{(p-1)^2}{p(p-2)})\prod_{p|m}(\frac{p}{p-1})\\
&=& C_0^{-1}\frac{m}{\phi(m)}.
\end{eqnarray*}
Moreover we have $m/\phi(m)\le e^{\gamma}\log x$ for $x\ge 9$, as shown
by Liu, Liu and Wang [13; (3.9)].  It then follows that
\begin{eqnarray*}
C_2&=&\int_{1}^{\infty}s(x)\frac{dx}{x^2}\\
&=&\int_{1}^{M}s(x)\frac{dx}{x^2}+\int_{M}^{\infty}s(x)\frac{dx}{x^2}\\
&\le&\sum_{\ep(d)\le M}\int_{\ep(d)}^{M}k(d)\frac{dx}{x^2}+
C_0^{-1}e^{\gamma}\int_{M}^{\infty}\log x\frac{dx}{x^2}\\
&\le &\sum_{\ep(d)<M}k(d)(\frac{1}{\ep(d)}-\frac{1}{M})+
2.744(\frac{1+\log M}{M})
\end{eqnarray*}
for any integer $M\ge 9$.

We now set
\[\sum_{\ep(d)=e}k(d)=\kappa(e)\]
so that
\[\sum_{e|d}\kappa(e)=\sum_{\ep(e)|d}k(e).\]
However $\ep(e)|d$ if and only if $e|2^d-1$.  Thus
\[\sum_{e|d}\kappa(e)=\sum_{e|2^d-1}k(e)=h(2^d-1).\]
We therefore deduce that
\[\kappa(e)=\sum_{d|e}\mu(e/d)h(2^d-1).\]
This enables us to compute
\[\sum_{\ep(d)<M}k(d)(\frac{1}{\ep(d)}-\frac{1}{M})=
\sum_{m<M}\kappa(m)(\frac{1}{m}-\frac{1}{M})\]
by using information on the prime factorization of $2^d-1$ for $d<M$. 
In particular, taking $M=20$ we find that
\[\sum_{m<20}\kappa(m)(\frac{1}{m}-\frac{1}{20})=1.6659\ldots,\]
and hence that
\begin{equation}
C_2\le\sum_{m<20}\kappa(m)(\frac{1}{m}-\frac{1}{20})+
2.744(\frac{1+\log 20}{20})=2.2141\ldots
\end{equation}
For comparison with this upper bound for $C_2$ we note that
\[C_2\ge\sum_{d\le 10000}k(d)/\ep(d)=1.9326\ldots\]
This latter figure is probably closer to the true value, but the
discrepancy is small enough for our purposes.

From (32), (40) and (41) we calculate that
\[(C_1-2)C_2+C_0^{-1}\log 2\le 13.967,\]
so that Lemma 9 yields the following bound.
\begin{lemma}
We have
\[J(\mathfrak{m})\le \{13.968+o(1)\}C_0\frac{N}{\log^2 2}.\]
\end{lemma}

\section{Completion of the Proof}

Let $R(N)$ denote the number of representations of $N$ as a sum of two
primes and $K$ powers of $2$ in the ranges under consideration, so that
\[R(N)=\int_0^1S(\alpha)^2 T(\alpha)^K e(-\alpha N)d\alpha.\]

To estimate the minor arc contribution to $R(N)$ we first bound $S(\alpha)$.
According to Theorem 3.1 of Vaughan [17] we have
\[\sum_{p\le x}e(\alpha p)\log p\ll
(\log x)^4\{xq^{-1/2}+x^{4/5}+x^{1/2}q^{1/2}\}\]
if $|\alpha-a/q|\le q^{-2}$ with $(a,q)=1$.  Thus if
$\alpha\in\mathfrak{m}$ we may take $P\ll q\ll N/P$ to deduce that
\[S(\alpha)\ll (\log N)^3\{N^{4/5}+NP^{-1/2}\}.\]
Taking $P=N^{45/154-4\vp}$, we obtain
\[S(\alpha)\ll N^{263/308+3\vp}.\]

If one assumes the Generalized Riemann
Hypothesis, we may apply Lemma 12 of Baker and Harman [1], which 
implies that
\[\sum_{n\le x}\Lambda(n)e((\frac{a}{q}+\beta) n)\ll 
(\log x)^2\{q^{-1}\min(x,|\beta|^{-1})+x^{1/2}q^{1/2}+x(q|\beta|)^{1/2}\}\]
when $|\beta|\le x^{-1/2}$.  It follows by partial summation that
\[S(\frac{a}{q}+\beta)\ll (\log N)\{q^{-1}\min(N,|\beta|^{-1})
+N^{1/2}q^{1/2}+N(q|\beta|)^{1/2}\}\]
for $|\beta|\le N^{-1/2}$.   According to Dirichlet's Approximation
Theorem, 
we can find $a$ and $q$ with
\[|\alpha-\frac{a}{q}|\le\frac{1}{qN^{1/2}},\;\;\;q\le N^{1/2}.\]
Thus
\[S(\alpha)\ll (\log N)N^{3/4}\]
unless $q\le N^{1/4}$ and $|\alpha-a/q|\le q^{-1}N^{-3/4}$.  Since
$\alpha\in\mathfrak{m}$ and
$P=N^{45/154-4\vp}\ge N^{1/4}$, these latter conditions cannot hold.

We therefore conclude that
\[S(\alpha)\ll N^{\theta+o(1)}\]
for $\alpha\in\mathfrak{m}$, where we take $\theta=263/308$ in general,
and $\theta=3/4$ under the Generalized Riemann Hypothesis.

We now have
\begin{eqnarray*}
\int_{{\mathfrak m}\cap\cl{A}_{\lambda}}S(\alpha)^2 T(\alpha)^K 
e(-\alpha N)d\alpha&\ll&{\rm
meas}(\cl{A}_{\lambda})N^{2\theta+o(1)}L^K\\
&\ll& N^{-E(\lambda)+2\theta+o(1)}\\
&\ll&N,
\end{eqnarray*}
providing that $E(\lambda)>2\theta-1$.  Thus, according to Lemma 1, we
may take $\lambda=0.863665$ unconditionally, and $\lambda=0.722428$
under the Generalized Riemann Hypothesis.

It remains to consider the set ${\mathfrak
m}\setminus\cl{A}_{\lambda}$.  Here we have
\begin{eqnarray*}
|\int_{{\mathfrak m}\setminus\cl{A}_{\lambda}}S(\alpha)^2 T(\alpha)^K 
e(-\alpha N)d\alpha|&\le& (\lambda L)^{K-2}
\int_{{\mathfrak m}}|S(\alpha)T(\alpha)|^2 d\alpha\\
&\le& (\lambda L)^{K-2} 13.968\frac{C_0}{\log^2 2} N.
\end{eqnarray*}

Finally we compare this with the estimate for the major arc integral,
given by Lemma 7, and conclude that
\[\int_0^1 S(\alpha)^2 T(\alpha)^K e(-\alpha N)d\alpha >0\]
providing that $N$ is large enough, $\varpi$ is small enough, and
\[
13.968\lambda^{K-2} < 2.7895.
\]
When $\lambda=0.863665$ this is satisfied for $K>12.991$, so that $K=13$ is
admissible.  Similarly, when $\lambda=0.722428$ one can take any $K>6.995$,
so that $K=7$ is admissible.  This completes the proof of our theorems,
subject to Lemma 1.

\section{Proof of Lemma 1}

In this section we shall prove Lemma 1.  We shall again use
$\vp$ to denote a small positive constant.  We shall allow the
constants implied by the $O(\ldots)$ and $\ll$ notations to depend on
$\vp$, although sometimes we shall mention the dependence explicitly
for emphasis.  As mentioned in the introduction, the
method we shall adopt was suggested to us by Professor Keith Ball, and 
is based on the martingale method for proving exponential inequalities
in probability theory.

It is convenient to work with 
\[T_L(\alpha)=T(\alpha/2)=\sum_{0\le n\le L-1}e(\alpha 2^n)\]
in place of $T(\alpha)$.  Clearly we have
\[{\rm meas}\{\alpha\in[0,1]: |T_L(\alpha)|\ge\lambda L\}=
{\rm meas}(\cl{A}_{\lambda}).\]
Let $M=1+[2\pi/\vp]$ and suppose that $|T_L(\alpha)|\ge\lambda L$ with
$\arg(T_L(\alpha))=\phi$.  Write $m=[M\phi/2\pi]$ and
$\rho_m=e(-m/M)$.  Then
\[|e^{-i\phi}-\rho_m|\le |\phi-\frac{2\pi m}{M}|\le
\frac{2\pi}{M}\le\vp,\]
whence
\begin{eqnarray*}
{\rm Re}(\rho_m T_L(\alpha))&\ge& {\rm Re}(e^{-i\phi}T_L(\alpha))-
\vp|T_L(\alpha)|\\
&=&(1-\vp)|T_L(\alpha)|\\
&\ge&(1-\vp)\lambda L.
\end{eqnarray*}
It follows that
\begin{eqnarray*}
\lefteqn{{\rm meas}\{\alpha\in[0,1]: |T_L(\alpha)|\ge\lambda L\}}
\hspace{2cm}\\
&\le&
\sum_{m=0}^{M-1}{\rm meas}\{\alpha\in[0,1]: {\rm Re}(\rho_m T_L(\alpha))
\ge(1-\vp)\lambda L\}\\
&\ll_{\vp}&\sup_{|\rho|=1}{\rm meas}\{\alpha\in[0,1]: 
{\rm Re}(\rho T_L(\alpha))\ge(1-\vp)\lambda L\}.
\end{eqnarray*}

We now set
\[S(\xi,\rho,L)=\int_{0}^{1}\exp\{\xi{\rm Re}(\rho
T_L(\alpha))\}d\alpha,\]
for an arbitrary real $\xi>0$, whence
\[S(\xi,\rho,L)\ge \exp\{\xi(1-\vp)\lambda L\}
{\rm meas}\{\alpha\in[0,1]: {\rm Re}(\rho
T_L(\alpha))\ge(1-\vp)\lambda L\}.\]
It therefore follows that
\begin{equation}
{\rm meas}(\cl{A}_{\lambda})\ll \exp\{-\xi(1-\vp)\lambda L\}
\sup_{|\rho|=1}S(\xi,\rho,L).
\end{equation}

For any integer $h$, we have
$T_L(\alpha)=T_{L-h}(2^h\alpha)+T_h(\alpha)$.  Moreover, for any function $f$
we have 
\[\int_0^1 f(\alpha)d\alpha = \frac{1}{2^h}\int_0^{1}
\sum_{r=0}^{2^h-1}f(\frac{\beta}{2^h}+\frac{r}{2^h})d\beta.\]
It therefore follows that
\[S(\xi,\rho,L)=\frac{1}{2^h}\int_0^{1}\sum_{r=0}^{2^h-1}
\exp\{\xi{\rm Re}(\rho T_{L-h}(\beta+r))\}
\exp\{\xi{\rm Re}(\rho T_{h}(\frac{\beta+r}{2^h}))\}d\beta.\]
Since $T(\alpha)$ has period $1$ this becomes
\[\int_0^{1}\exp\{\xi{\rm Re}(\rho T_{L-h}(\beta))\}\frac{1}{2^h}
\sum_{r=0}^{2^h-1}\exp\{\xi{\rm Re}(\rho T_{h}(\frac{\beta+r}{2^h}))\}
d\beta.\]
If we now set
\begin{equation}
F(\xi,h)=\sup_{\beta\in[0, 1],\, |\rho|=1}\frac{1}{2^h}
\sum_{r=0}^{2^h-1}\exp\{\xi{\rm Re}(\rho T_h(\frac{\beta +
r}{2^h}))\}
\end{equation}
we deduce that
\[S(\xi,\rho,L)\le S(\xi,\rho,L-h)F(\xi,h).\]
Using this inductively we find that
\[S(\xi,\rho,L)\le S(\xi,\rho,L-nh)F(\xi,h)^n,\]
and taking $n=[L/h]$ we deduce that
\[S(\xi,\rho,L)\ll_{\xi,h} F(\xi,h)^n\ll_{\xi,h} F(\xi,h)^{L/h}.\]
When we combine this with (42) we deduce that
\[{\rm meas}(\cl{A}_{\lambda})\ll_{\xi,h,\vp} \exp\{-\xi(1-\vp)\lambda L\}
F(\xi,h)^{L/h}.\]
It follows that we may take
\[E(\lambda)=\frac{\xi\lambda}{\log 2}-\frac{\log F(\xi,h)}{h\log 2}-
\frac{\vp}{\log 2}\]
for any $h\in\N$, any $\xi>0$ and any $\vp>0$.  

We proceed to show that the supremum in (43) occurs at $\beta=0$ and
$\rho=1$, whence
\begin{equation}
F(\xi,h)=\frac{1}{2^h}\sum_{r=0}^{2^h-1}
\exp\{\xi{\rm Re}(T_h(\frac{r}{2^h}))\}.
\end{equation}
Since
\[{\rm Re}(\rho T_h(\frac{\beta+r}{2^h}))=\frac{1}{2}
\{\rho T_h(\frac{\beta+r}{2^h})+\overline{\rho}\,T_h(\frac{-\beta-r}{2^h})\},\]
we find that
\begin{eqnarray*}
\lefteqn{\sum_{r=0}^{2^h-1}\exp\{\xi{\rm Re}(\rho T_h(\frac{\beta+r}{2^h}))\}}
\hspace{2cm}\\
& = & \sum_{n=0}^\infty\frac{1}{2^n\cdot n!}\sum_{r=0}^{2^h-1}
\xi^n \left(\rho T_h(\frac{\beta + r}{2^h}) + \overline{\rho}\,
T_h(\frac{-\beta-r}{2^h})\right)^n.
\end{eqnarray*}
However
\[\sum_{r=0}^{2^h-1}
\left(\rho T_h(\frac{\beta + r}{2^h}) + \overline{\rho}\,
T_h(\frac{-\beta-r}{2^h})\right)^n
  =  \sum_{m=0}^n \left(\begin{array}{c}n\\m\end{array}\right) 
\rho^{2m-n} S(n,m,h,\beta),\]
where
\begin{equation}
S(n,m,h,\beta)=\sum_{r=0}^{2^h-1}
T_h(\frac{\beta+r}{2^h})^m  T_h(\frac{-\beta-r}{2^h})^{n-m}.
\end{equation}
It follows that
\begin{equation}
F(\xi,h)\le \frac{1}{2^h}\sup_{\beta\in[0,1]}
\sum_{n=0}^\infty\frac{1}{2^n\cdot n!}
\xi^n \sum_{m=0}^n {n\choose m} |S(n,m,h,\beta)|.
\end{equation}

We now expand the powers of $T_h$ occurring in (45), 
and perform the summation over $r$. 
We then see that $S(n,m,h,\beta)$ is a sum of terms
\[2^h\exp\{\beta(2^{a_1}+\ldots+2^{a_m}-2^{b_1}-\ldots-2^{b_{n-m}})\},\]
over integer values $a_i,b_j$ between 0 and $h-1$, subject to the
condition
\[2^{a_1}+\ldots+2^{a_m}\equiv 2^{b_1}+\ldots+2^{b_{n-m}}\pmod{2^h}.\]
It is now apparent that $|S(n,m,h,\beta)|\le S(n,m,h,0)$, whence (46)
yields
\begin{eqnarray*}
F(\xi,h)&\le &\frac{1}{2^h}\sum_{n=0}^\infty\frac{1}{2^n\cdot n!}
\xi^n \sum_{m=0}^n {n\choose m} S(n,m,h,0)\\
&=&\frac{1}{2^h}\sum_{r=0}^{2^h-1}
\exp\{\xi{\rm Re}(T_h(\frac{r}{2^h}))\},
\end{eqnarray*}
The assertion (44) now follows.

Hence it remains to compute $F(\xi,h)$ using (44)
and optimize for $\xi$ in (42).  We have carried out the computations 
for $h=16$. Comparing the results
for this value with the outcome for smaller values of $h$, it appears
that the potential improvements obtainable by choosing $h$ larger 
than 16 are only small.  After taking suitable care over rounding
errors we find that we may take $\xi=1.181$ to get
\[E(0.863665)>\frac{109}{154}+10^{-8}\]
and $\xi=0.905$ to get
\[E(0.722428)>\frac{1}{2}+10^{-8}.\]
Using Mathematica 4.1 on a PC, computing the values
$T_{16}(r/2^{16})$ for the integers
$0\leq r\leq 2^{16}-1$ took about 7 minutes, and
summing these values up to obtain $F(\xi, h)$ took 24 seconds for each
of the two values of $\xi$.

\bigskip
\bigskip

D.R. Heath-Brown

Mathematical Institute,

24-29, St.Giles',

Oxford OX1 3LB,

ENGLAND
\bigskip

rhb@maths.ox.ac.uk
\bigskip
\bigskip
\bigskip

J.-C. Puchta

Mathematical Institute,

24-29, St.Giles',

Oxford OX1 3LB,

ENGLAND
\bigskip

puchta@maths.ox.ac.uk

\end{document}